\documentclass[11pt,fleqn,centertags]{article}
\usepackage{mathptmx}
\usepackage{amsmath}
\usepackage{amssymb}
\usepackage{amsthm}
\usepackage[mathscr]{eucal}
\usepackage[all,tips]{xy}
\usepackage{hyperref}

\newtheorem{Theorem}{Theorem}
\newtheorem{rmk}{Remark}

\theoremstyle{definition}

\DeclareMathOperator{\PSL}{PSL}

\def\Hy{\mathbb{H}}
\def\PP{\mathcal{P}}

\def\C{{\mathbb C}}

\def\D{\mathcal{D}}

\begin{document}
\title{\textbf{On fields of definition of arithmetic Kleinian reflection groups}}
\author{Mikhail Belolipetsky\thanks{Belolipetsky partially supported by EPSRC grant EP/F022662/1}}
\date{\empty}
\maketitle

\begin{abstract}
We show that degrees of the real fields of definition of arithmetic Kleinian reflection
groups are bounded by $35$.
\end{abstract}

\section{Introduction}

In two recent articles \cite{N1}, \cite{N2}, Nikulin gave explicit upper bounds
for the degrees of the fields of definition of arithmetic hyperbolic reflection
groups in dimensions $n\ge 4$ and $n = 2$. He pointed out that only the case $n=3$
remained open. Later, in \cite{N3}, Nikulin extended his method to this case as
well and thus completed in general the solution of the problem. Still the explicit
bounds obtained in \cite{N1,N2,N3} are far from being sharp.

In this note we explore an alternative approach to the problem. As a result we
show that fields of definition of arithmetic Kleinian reflection groups have
degrees less than $70$. This implies that real fields of definition of
corresponding orthogonal groups have degrees bounded by $35$ (compare with
$10000$-bound in \cite{N3}, which was improved to $909$ in the latest version of
the preprint). We also give general bounds for discriminants
of the fields of definition and good upper bounds for the discriminants which
correspond to non-cocompact Kleinian groups.

A theorem of Nikulin \cite[Th.~4.8]{N1} states that the degrees of the real
fields of definition in all dimensions are bounded by maximum of $56$, and the
degrees in dimensions $2$ and $3$. In dimension $2$, following the previous
work \cite{LMR}, Nikulin gave an upper bound $44$ for the degrees (see
\cite[Sec.~4.5]{N1}). Combining this with the result of the present paper we
get a universal upper bound $56$ for all dimensions, which, in particular,
improves on the bounds for the dimensions $4$ and $5$ obtained by Nikulin in
\cite{N2}.

Our method is based on the work of Agol \cite{A1} combined with Borel's volume
formula \cite{B1} and some number-theoretic results of Chinburg and Friedman
\cite{C1,CF}.

\section{Preliminaries}

Discrete subgroups of $\PSL(2,\C)$ are called {\em Kleinian groups}. As
$\PSL(2,\C)$ is isomorphic to the group of orientation preserving isometries of
the hyperbolic $3$-space $\Hy^3$, Kleinian groups act isometrically on $\Hy^3$.
If a Kleinian group $\Gamma$ acts as an orientation preserving subgroup of a
discrete group generated by reflections in hyperbolic hyperplanes, it is called
a {\em reflection group}. In this case the volume of the hyperbolic polyhedron
$\PP$ bounded by the reflection hyperplanes is half of the covolume of
$\Gamma$. If $\Gamma$ is an arithmetic subgroup of $\PSL(2,\C)$ (we refer to
\cite[Sec.~8]{MR} for the definition and basic properties of arithmetic Kleinian
groups), its covolume can be estimated using its arithmetic invariants. A volume
formula of Borel \cite{B1} allows us to write down the estimate in an explicit
form. On the other hand, the covolume of $\Gamma$ can be estimated from the
geometry of the polyhedron $\PP$. One of the main ideas of $\cite{A1}$ is that
the interplay between these two estimates leads to the finiteness result for
the number of conjugacy classes of arithmetic Kleinian reflection groups.
Our purpose is to make this relation explicit and then apply it to get
quantitative bounds for the finiteness theorem.

\section{Results}
\begin{Theorem}
Let $k$ be a field of definition of an arithmetic Kleinian reflection group.
Then its degree $n_k \le 70$ and the absolute value of discriminant $\D_k < 4.4\times10^{273}$.
Moreover, if the group is non-cocompact then $n_k = 2$ and $\D_k\le 9240$.
\end{Theorem}

\begin{rmk} It is well-known due to Vinberg that arithmetic reflection groups
are defined by quadratic forms \cite[Lemma 7]{Vinb1}. Corresponding
forms are defined over totally real fields, and the fields of definition
of arithmetic Kleinian reflection groups are quadratic extensions of these
fields (see \cite[Sec.~10.2]{MR}). Therefore, the degrees of the fields
of definition of arithmetic Kleinian reflection groups are necessarily even
and corresponding real fields of definition of orthogonal groups have
degrees twice less. Theorem~1 implies that the degrees of the real
fields of definition are bounded by $70/2 = 35$.
\end{rmk}

\begin{rmk}
The non-cocompact arithmetic Kleinian groups are also known as Bianchi groups
and have a long history in mathematical literature. While our general bound for
discriminants of the fields of definition is more an existence bound, the bound
for the Bianchi groups is much better (compare with $1.02\times10^8$-bound in
\cite[Sec.~7]{A1}). Moreover, in the proof of the theorem we obtained a list of
$330$ imaginary quadratic fields which contains all possible candidates for the
fields of definition of Bianchi reflection groups. The key ingredients for this
enumeration are the Gauss' theorem on the $2$-class number of a quadratic field
and Vinberg's result on class groups of fields of definition of Bianchi
reflection groups (see part 3 of the proof for details). No proper analogue of
any of these two results is currently known in a more general setting.
\end{rmk}

\section{Proof of Theorem 1}

\subsection{}

Let $\Gamma$ be a maximal arithmetic Kleinian reflection group and let
$\mu(\Gamma)$ denote its covolume. By the proof of Theorem~6.1 in \cite{A1},
we have
\begin{equation}
\mu(\Gamma) \le 2\times 64\pi^2.
\end{equation}
We now recall a non-trivial corollary of Borel's volume formula which was
obtained in \cite[Lemma 4.3]{CF}. It implies
\begin{equation}
\mu(\Gamma) > 0.69\exp\left(0.37n_k - \frac{19.08}{h(k,2,B)}\right),
\end{equation}
where $n_k$ is the degree of the field of definition of $\Gamma$ and $h(k,2,B)$
is the order of a certain subgroup of the $2$-class group of $k$ (see
\cite[Sec.~2]{CF} for a precise definition of $h(k,2,B)$).

From (1) and (2) we obtain
\begin{align*}
n_k & < \frac{\log(128\pi^2/0.69) + 19.08/h(k,2,B)}{0.37}\\
    & \le \frac{\log(128\pi^2/0.69) + 19.08}{0.37}\\
    & < 71.88.
\end{align*}

It was already pointed out in Remark~1 that the degree of $k$ is necessarily
an even integer, which implies $n_k \le 70$. As any arithmetic reflection
group is contained in some maximal arithmetic reflection group it follows that
the same bound holds for the fields of definition of arbitrary arithmetic
Kleinian reflection groups.

\subsection{}

Using the Brauer-Siegel theorem and Zimmert's bound for the regulator we can show
that the class number
\begin{equation}
h_k \le 10^2 \left(\frac{\pi}{12}\right)^{n_k} \D_k
\end{equation}
(see \cite[Sec.~3]{CF} for more details).

The bound (3) together with the upper bound for $n_k$ gives
\begin{equation*}
h(k,2,B) \le h_k \le c_1\D_k,\quad c_1 = 10^2 \left(\frac{\pi}{12}\right)^{70}.
\end{equation*}

Now the volume formula \cite{B1} implies (here again we use the notation from \cite{CF}):
\begin{align}
\mu(\Gamma) &\ge
\frac{\D_k^{3/2} \zeta_k(2) 2^{t-t'}}{2^{2(n_k-2)+4+t-1}\pi^{2(n_k-2)+2}[K(B):k]}
\prod_{\substack{v\in R_f \\ 2\neq Nv}}\frac{(Nv-1)}{2} \\
&\ge \frac{\D_k^{3/2}}{2^{t'}2^{2n_k-1}\pi^{2n_k-2}2^{2(n_k-2)} h(k,2,B)} \nonumber\\
&\ge \frac{\D_k^{1/2}}{c_2 c_1},\quad c_2 = 2^{70}2^{4\times 70-5}\pi^{2\times 70-2} \nonumber
\end{align}

Combined with (1) this gives
\begin{equation*}
\D_k \le (128\pi^2 c_1 c_2)^2 < 4.4\times 10^{273}.
\end{equation*}

\subsection{}

Let now $\Gamma$ be a non-cocompact arithmetic Kleinian group. It follows that
$k$ is necessarily imaginary quadratic field and the set of places of $k$ at which
the quaternion division algebra associated to $\Gamma$ is ramified is empty (see
\cite[Th.~8.2.3]{MR}). In the notation of \cite{CF} the latter implies $t'=0$.
Thus the volume estimate (4) applied to this case gives
\begin{equation}
\mu(\Gamma) \ge \frac{\D_k^{3/2} \zeta_k(2)}{2^3\pi^2 h(k,2,B)}.
\end{equation}
We have $h(k,2,B) \le h_2(k)$, where $h_2(k)$ is the $2$-class number of $k$.
By Gauss' theorem, $h_2(k) \le 2^{t_k-1}$, where $t_k$ is the number of distinct
prime divisors of the discriminant $\D_k$. Well-known number-theoretic
estimates imply $t_k \le 1.5\log\D_k/\log\log\D_k$.
Therefore,
\begin{align*}
\mu(\Gamma) &\ge \frac{\D_k^{3/2} \zeta_k(2)}{2^{1.5\log\D_k/\log\log\D_k - 1}8\pi^2} \\
            &\ge \frac{\D_k^{3/2}}{2^{1.5\log\D_k/\log\log\D_k}16\pi^2}.
\end{align*}
If $\D_k \ge 10^5$, then $\mu(\Gamma) > 1492.9$ and inequality (1) fails. Considering
the remaining fields, a simple program for GP PARI calculator allows us to compute the
lower bound (5) for $\mu(\Gamma)$ using precise values of $h_2(k)$ and $\zeta_k(2)$.
This way we obtain that there are in total $882$ fields which satisfy the criteria,
and that the largest $\D_k = 9240$. The list of the admissible fields can be further
improved using a theorem of Vinberg \cite{Vinb2} which implies that the class numbers
of the fields under consideration should be powers of $2$. As a result we obtain that
there are only $330$ such fields.
\medskip

\noindent {\it Acknowledgement.} I would like to thank Ian Agol for discussions
and the referee for helpful comments.

\noindent Department of Mathematical Sciences, Durham University, Durham DH1 3LE, UK \\
Sobolev Institute of Mathematics, Koptyuga 4, 630090 Novosibirsk, RUSSIA \\
email: {\tt mikhail.belolipetsky@durham.ac.uk}

\end{document}